# Analytical and numerical investigation of radiative heat transfer in semitransparent Media

Yaochuang Han

**Abstract**

This paper is devoted to deal with some mathematical and numerical aspects of the radiative integral transfer equations. First, the properties of the raidative integral operators are analyzed. Based on these results, the existence and uniqueness of solution to the radiative integral system is proved by the reversibility of operator matrix. Besides, the convergence analysis of an iterative scheme is also carried out. Then, boundary element method based on the collocation scheme is used to discretize the radiative integral system. Arbitrary enclosure geometry is considered. For the non-convex geometries, an element-subdivision algorithm is developed to handle with the integrals contained the visibility factor. Finally, examples are presented to verify the effectiveness and accuracy of our method.

Key words: Radiative heat transfer, Iteration method, Boundary element method, Shadow detection

**Introduction**

It is well known that radiative heat transfer dominates the energy transport in many engineering applications, such as astrophysics, nuclear reactors, reentry of space vehicles, combustion in gas turbine combustion chambers, etc. As a result, research on the radiative transfer problem has always been a very active and important area.

In the last decades, there are two kinds of models depicting the radiative heat transfer problem, which are the integro-differential radiative transfer equation (IDRTE) and the radiative integral transfer equations (RITEs) respectively. For the former, some comprehensive summarizations about mathematical theory and numerical methods can be found in [1] and [2, 3] respectively. Because of inherent severe nonlinearity and high dimensionality, the fast iterative algorithm is critical in simulation, of which [4] provides valuable reference. The RITEs are obtained by integrating the IDRTE over the entire solid angles. As a consequence, the RITEs only contain the spatial variables. This technique reduces the dimensions at a cost of singular kernels and dense matrices. Much work involved in numerical simulation for RITEs has been investigated. Crosbie [5, 6] deduced multi-dimensional RITEs in rectangular enclosures. Thynell [7] formulated the RITEs in two-dimensional cylindrical enclosures, with absorbing, emitting, and linear-anisotropically scattering media. A set of integral equations for thermal radiative transmission in anisotropic scattering media is formulated and discretized by the product-integration in [8]. Based on an alternative form of formulation in [8], Sun [9] developed a modified boundary element method to simulate the thermal radiation problems. Li et al [10] adopted the Galerkin boundary element method to discretize the thermal radiation problem without scattering. Altac [11, 12] used a method of "subtraction of singularity" to calculate the surface integrals and volume integrals in RITEs. Compared with numerical simulation, few attentions are focused on theory analysis for RITEs. To our knowledge, the theoretical researches about RITEs are confined to the *radiosity equation*, which is the simplest version of RITEs. The radiosity equation has been deeply investigated in [13-16]. The collocation method is proposed and analyzed in [13], where the unique solvability of radiosity equation is considered by properties of the radiosity integral operator. For a discussion of the numerical integration needed in implementing collocation methods, see [14]. The regularity of the solution on polyhedral surfaces is studied in [15]. In [16], the authors consider the rigorous convergence and error analysis of the Galerkin boundary element method in convex and non-convex enclosure geometries. Meanwhile, the coupling of heat radiation to other heat transfer mechanisms has also received some attentions, such as [17, 18]. Laitinen [17] gave the existence of a weak solution to the heat equation for a non-convex body with Stefan-Boltzmann radiation condition on the surface. Following that, they analyzed some mathematical results for coupled conduction-radiation heat transfer in semitransparent materials in [18]. The mathematical respects for RITEs in scattering media seem to be ignored by the mathematical community. It is no doubt that the theoretical research for RITEs is of importance for further coupled energy transfer. So, this paper is devoted to make some contributions to this point.

Additionally, for the non-convex geometries, the detection of the shadow zone seriously affects the accuracy of the final results, and it is very time-consuming. As a result, it is a key point for the heat radiation calculation. Research on this topic is extensive. The monographs about this topic can be refered to [19, 20]. Reference [21] gives a summary and comparison of existing algorithms used in heat radiation calculations. The analogous researches can be found [22, 23]. In [22], the authors only considered 2D problem. Additionally, the subdivision of the current element relies on the visibility of the all Gaussian nodes in the element. The discrete distribution of Gaussian nodes brings in unnecessary errors in this step. The algorithm in [23] is based on the Galerkin boundary element method. As is said in [27], '*The Galerkin technique requires double integration over surface and volume. Thus to calculate the entries of the matrices, standard zoning methods require expensive integration in four, five or six dimensions.*' The present algorithm is developed based on the idea of [22]. In order to improve the accuracy, an idea of element-subdivision is introduced to handle with the visibility factor in the integral operators.

The outline of this paper is as follows. In the second section, the formulations and some physical background of RITEs are presented. Follow that, the properties of four radiative integral operators are analyzed. Then the existence and uniqueness of the radiative integral formulae is proved based on the reversibility of operator matrix in section 3. In section 4, the convergence of an iterative scheme is proved based on the properties of the radiative integral operators with appropriate physical parameters. Also the boundary element method based on the collocation scheme is briefly described. In section 5, a new high-precision visibility algorithm is presented in detail. In section 6, two numerical examples on convex and non-convex domains are tested by using the aforementioned iterative scheme and boundary element method. Finally, the paper ends with some concluding remarks.

**2. Operator expression of RITEs**

Consider radiation transport in the bounded domain $V \subset R^3$ with boundary $S$. The boundary is assumed to be diffuse and gray. That is to say, the emissivity and absorptivity of the surface are independent of direction and wavelength of the radiation. Additionally, the medium in enclosure is homogeneous and isotropic.

Under the above simplifications, RITEs in an absorbing, emitting, and scattering medium can be written as follows [11]

$$q(\mathbf{p})+\varepsilon E_b(\mathbf{p})=\varepsilon\int_S\left[E_b(\mathbf{r})+\frac{1-\varepsilon}{\varepsilon}q(\mathbf{r})\right]\exp[-\beta|\mathbf{p}-\mathbf{r}|]\frac{\cos\phi_p\cos\phi_r}{\pi|\mathbf{p}-\mathbf{r}|^2}\chi(\mathbf{p},\mathbf{r})dS(\mathbf{r})+$$
$$\varepsilon\sigma_a\int_V I_b(\mathbf{r})\exp(-\beta|\mathbf{p}-\mathbf{r}|)\frac{\cos\phi_p}{|\mathbf{p}-\mathbf{r}|^2}\chi(\mathbf{p},\mathbf{r})dV(\mathbf{r})+$$
$$\frac{\varepsilon\sigma_s}{4\pi}\int_V G(\mathbf{r})\exp(-\beta|\mathbf{p}-\mathbf{r}|)\frac{\cos\phi_p}{|\mathbf{p}-\mathbf{r}|^2}\chi(\mathbf{p},\mathbf{r})dV(\mathbf{r}),\qquad \mathbf{p}\in S, \quad (1)$$

$$G(\mathbf{p})=\int_S\frac{1}{\pi}\left[E_b(\mathbf{r})+\frac{1-\varepsilon}{\varepsilon}q(\mathbf{r})\right]\exp(-\beta|\mathbf{p}-\mathbf{r}|)\frac{\cos\phi_r}{|\mathbf{p}-\mathbf{r}|^2}\chi(\mathbf{p},\mathbf{r})dS(\mathbf{r})$$
$$+\sigma_a\int_V I_b(\mathbf{r})\exp(-\beta|\mathbf{p}-\mathbf{r}|)\frac{1}{|\mathbf{p}-\mathbf{r}|^2}\chi(\mathbf{p},\mathbf{r})dV(\mathbf{r})$$
$$+\frac{\sigma_s}{4\pi}\int_V G(\mathbf{r})\exp(-\beta|\mathbf{p}-\mathbf{r}|)\frac{1}{|\mathbf{p}-\mathbf{r}|^2}\chi(\mathbf{p},\mathbf{r})dV(\mathbf{r}),\qquad \mathbf{p}\in V, \quad (2)$$

where $q(\mathbf{p}), E_b(\mathbf{p}), I_b(\mathbf{p}), G(\mathbf{p})$ denote the radiative flux, blockbody emissive power, blockbody intensity of radiation and the incident energy at $\mathbf{p}$, respectively. $\varepsilon$ represents the emissivity of boundary. $\sigma_a, \sigma_s, \beta$ are the absorptivity, scattering coefficient, extinction coefficient with $\beta=\sigma_a+\sigma_s$ within the medium. $\phi$ stands for the angle between the incoming ray and the outward normal to the surface.

The Boolean function $\chi$ in Eqs. (1), (2), named as the shadow zone function, is defined as
$$\chi(\mathbf{p},\mathbf{r})=\begin{cases}1, & \text{if } \mathbf{r} \text{ can be seen by } \mathbf{p},\\ 0, & \text{otherwise.}\end{cases} \quad (3)$$
Here the statement "$\mathbf{r}$ can be seen by $\mathbf{p}$" means that there is no opaque material between $\mathbf{r}$ and $\mathbf{p}$, (i.e. $\overline{\mathbf{rp}}\cap S=\varnothing$).

Blockbody emissive power can be computed from the *Stefan-Boltzmann law*
$$E_b(\mathbf{p})=\sigma T^4(\mathbf{p}),\qquad I_b(\mathbf{p})=\frac{\sigma T^4(\mathbf{p})}{\pi},$$
where $\sigma$ denote the Stefan-Boltzmann constant, $T(\mathbf{p})$ is the temperature at $\mathbf{p}$.

So, once temperatures of both the medium and the bounding surface are known, Eqs. (1), (2) can be rewritten in operator form as
$$K\mathbf{u}=\mathbf{f} \quad (4)$$
where
$$K=\begin{pmatrix}I-K_1 & -K_2\\ -K_4 & I-K_3\end{pmatrix},\quad \mathbf{u}=\begin{pmatrix}q\\ G\end{pmatrix},\quad \mathbf{f}=\begin{pmatrix}f_1\\ f_2\end{pmatrix} \quad (5)$$

$$(K_1 q)(\boldsymbol{p}) = \frac{(1-\varepsilon)}{\pi} \int_S \exp\left[-\beta|\boldsymbol{p}-\boldsymbol{r}|\right] \frac{\cos\phi_p \cos\phi_r}{|\boldsymbol{p}-\boldsymbol{r}|^2} \chi(\boldsymbol{p},\boldsymbol{r}) q(\boldsymbol{r}) dS(\boldsymbol{r})$$

$$(K_2 G)(\boldsymbol{p}) = \frac{\varepsilon \sigma_s}{4\pi} \int_V \exp(-\beta|\boldsymbol{p}-\boldsymbol{r}|) \frac{\cos\phi_p}{|\boldsymbol{p}-\boldsymbol{r}|^2} \chi(\boldsymbol{p},\boldsymbol{r}) G(\boldsymbol{r}) dV(\boldsymbol{r})$$

$$(K_3 G)(\boldsymbol{p}) = \frac{\sigma_s}{4\pi} \int_V \exp(-\beta|\boldsymbol{p}-\boldsymbol{r}|) \frac{1}{|\boldsymbol{p}-\boldsymbol{r}|^2} \chi(\boldsymbol{p},\boldsymbol{r}) G(\boldsymbol{r}) dV(\boldsymbol{r})$$

$$(K_4 q)(\boldsymbol{p}) = \frac{1-\varepsilon}{\varepsilon \pi} \int_S \exp(-\beta|\boldsymbol{p}-\boldsymbol{r}|) \frac{\cos\phi_r}{|\boldsymbol{p}-\boldsymbol{r}|^2} \chi(\boldsymbol{p},\boldsymbol{r}) q(\boldsymbol{r}) dS(\boldsymbol{r})$$

(6)

$$f_1(\boldsymbol{p}) = \varepsilon \int_S E_b(\boldsymbol{r}) \exp\left[-\beta|\boldsymbol{p}-\boldsymbol{r}|\right] \frac{\cos\phi_p \cos\phi_r}{\pi|\boldsymbol{p}-\boldsymbol{r}|^2} \chi(\boldsymbol{p},\boldsymbol{r}) dS(\boldsymbol{r}) - \varepsilon E_b(\boldsymbol{p}) +$$

$$\varepsilon \sigma_a \int_V I_b(\boldsymbol{r}) \exp(-\beta|\boldsymbol{p}-\boldsymbol{r}|) \frac{\cos\phi_p}{|\boldsymbol{p}-\boldsymbol{r}|^2} \chi(\boldsymbol{p},\boldsymbol{r}) dV(\boldsymbol{r})$$

$$f_2(\boldsymbol{p}) = \int_S E_b(\boldsymbol{r}) \exp(-\beta|\boldsymbol{p}-\boldsymbol{r}|) \frac{\cos\phi_r}{\pi|\boldsymbol{p}-\boldsymbol{r}|^2} \chi(\boldsymbol{p},\boldsymbol{r}) dS(\boldsymbol{r}) +$$

$$\sigma_a \int_V I_b(\boldsymbol{r}) \exp(-\beta|\boldsymbol{p}-\boldsymbol{r}|) \frac{1}{|\boldsymbol{p}-\boldsymbol{r}|^2} \chi(\boldsymbol{p},\boldsymbol{r}) dV(\boldsymbol{r})$$

(7)

and $I$ denotes the identity operator.

The properties of each operator will be examined in the next section separately.

## 3. Some mathematical results of RITEs

The integral system Eq. (4) is a coupled Fredholm integral system of the second kind. As is said in [24], '*Existence and uniqueness of a solution to an operator equation can be equivalently expressed by the existence of the inverse operator.*' So a way to establish existence and uniqueness of the system Eq. (4) is to prove the existence of the inverse operator of $K$. Naturally, the properties of operators $K_1, K_2, K_3, K_4$ need to be clear, which are presented in the following.

Firstly, we extend the Lemma 2.1 in [13] to non-convex domain. Then we have the following Lemma.

**Lemma 1.** Suppose $V$ is a bounded domain of $R^3$ and has a Lipschitz boundary $S$. Let $\boldsymbol{p} \in S$, and let $S$ be smooth in an open neighborhood of $\boldsymbol{p}$. Besides, $\boldsymbol{p}$ can see all the other points in the closure $\overline{V}$. Then

$$\int_S \frac{\cos\phi_p \cos\phi_r}{|\boldsymbol{p}-\boldsymbol{r}|^2} dS(\boldsymbol{r}) = \pi. \tag{8}$$

Proof. Let $\delta$ be a sufficiently small number. Exclude an $\delta-$neighborhood of $\boldsymbol{p}$ from $V$,

and denote the remaining set by $V'$:

$$V' = V \setminus \{r \in V : |r - p| \leq \delta\}.$$

Let $S'$ denote the boundary of $V'$, and let $S_\delta$ denote the boundary of $V \setminus V'$. Then

$$\int_S \frac{\cos\phi_p \cos\phi_r}{|p-r|^2} dS(r) = \int_{S'} \frac{\cos\phi_p \cos\phi_r}{|p-r|^2} dS(r) + \int_{S_\delta} \frac{\cos\phi_p \cos\phi_r}{|p-r|^2} dS(r).$$

For a continuously differentiable vector function $v(r)$ defined $V'$, the divergence theorem says

$$\int_{S'} v(r) \cdot n_r dS(r) = -\int_{V'} \nabla \cdot v(r) dV(r).$$

We apply this with

$$v(r) = \frac{(p-r) \cdot n_p}{|r-p|^4}(r-p).$$

A straightforward computation shows

$$\nabla \cdot v(r) = 0, \quad r \in V',$$

And, therefore,

$$\int_{S'} \frac{\cos\phi_p \cos\phi_r}{|p-r|^2} dS(r) = \int_{S'} \frac{[(p-r)\cdot n_p][(r-p)\cdot n_r]}{|p-r|^4} dS(r) = 0.$$

Decompose $S_\delta$ into two parts:

$$S_\delta = T_\delta \cup W_\delta,$$

with

$$T_\delta = \{r \in S \,|\, |r-p| \leq \delta\},$$
$$W_\delta = \{r \in V \,|\, |r-p| = \delta\}.$$

Then

$$\int_{S_\delta} \frac{\cos\phi_p \cos\phi_r}{|p-r|^2} dS(r) = \int_{T_\delta} \frac{\cos\phi_p \cos\phi_r}{|p-r|^2} dS(r) + \int_{W_\delta} \frac{\cos\phi_p \cos\phi_r}{|p-r|^2} dS(r).$$

According to the assumption, that is, $p$ can see all the other points in the closure $\overline{V}$. Then the following inequalities are valid

$$0 \leq \phi_p, \phi_r \leq \frac{\pi}{2},$$

which is followed by

$$\frac{\cos\phi_p \cos\phi_r}{|p-r|^2} \geq 0.$$

Also, from [12] we have

$$\frac{\cos\phi_p \cos\phi_r}{|p-r|^2} \leq c, \quad p,r \in S, p \neq r,$$

with $c$ independent of $p$ and $r$.

Then

$$0 \leq \int_{T_\delta} \frac{\cos\phi_p \cos\phi_r}{|p-r|^2} dS(r)$$

$$\leq c \int_{T_\delta} 1 dS(r)$$

$$= O(\delta^2).$$

Thus this integral goes to zero as $\delta \to 0$.

For any $r \in W_\delta$,

$$n_r = \frac{r-p}{|r-p|}, \quad n_r \cdot \frac{r-p}{|r-p|} = 1.$$

So,

$$\int_{W_\delta} \frac{\cos\phi_p \cos\phi_r}{|p-r|^2} dS(r) = \int_{W_\delta} \frac{[(p-r)\cdot n_p][(r-p)\cdot n_r]}{|p-r|^4} dS(r)$$

$$= \int_{W_\delta} \frac{(p-r)\cdot n_p}{|p-r|^3} dS(r)$$

$$= \frac{1}{\delta^3} \int_{W_\delta} (p-r)\cdot n_p dS(r)$$

$$= \frac{1}{\delta^3} \int_0^{2\pi} \int_0^{\pi/2} \delta\cos\varphi \delta^2 \sin\varphi d\varphi d\theta$$

$$= \pi.$$

In conclusion, the proof is finished.

**Lemma 2**. Assume $S$ is a piecewise smooth surface of bounded domain $V$ in $R^3$. Let $p \in S$ be a point at which $S$ is smooth. Then

$$\int_S \exp[-\beta|p-r|] \frac{\cos\phi_p \cos\phi_r}{|p-r|^2} \chi(p,r) dS(r) \leq \pi, \tag{9}$$

$$\int_V \exp(-\beta|p-r|) \frac{\cos\phi_p}{|p-r|^2} \chi(p,r) dV(r) \leq \frac{\pi}{\beta}. \tag{10}$$

*Proof.* For any two points $p, r$, we have

$$\exp[-\beta|p-r|] \leq 1.$$

As a consequence, the following inequality is valid

$$\int_S \exp[-\beta|\pmb{p}-\pmb{r}|]\frac{\cos\phi_p \cos\phi_r}{|\pmb{p}-\pmb{r}|^2}\chi(\pmb{p},\pmb{r})dS(\pmb{r}) \leq \int_S \frac{\cos\phi_p \cos\phi_r}{|\pmb{p}-\pmb{r}|^2}\chi(\pmb{p},\pmb{r})dS(\pmb{r}).$$

For any $\pmb{p} \in S$, we always can find a bounded domain $V' \subset V$ with boundary $S'$, so that $\pmb{p}$ can see all the other point in the closure $\overline{V'}$.

So,

$$\int_S \frac{\cos\phi_p \cos\phi_r}{|\pmb{p}-\pmb{r}|^2}\chi(\pmb{p},\pmb{r})dS(\pmb{r}) = \int_{S'} \frac{\cos\phi_p \cos\phi_r}{|\pmb{p}-\pmb{r}|^2}dS(\pmb{r}).$$

Apply the preceding lemma 1, and we have

$$\int_{S'} \frac{\cos\phi_p \cos\phi_r}{|\pmb{p}-\pmb{r}|^2}dS(\pmb{r}) = \pi.$$

In conclusion, the inequality (6) holds.

As shown in [25],

$$dL_{rp}(\pmb{r}')\frac{dS(\pmb{r})\cos\phi_r}{|\pmb{r}-\pmb{p}|^2} = \frac{dV(\pmb{r}')}{|\pmb{r}'-\pmb{p}|^2}.$$

Adopting this equality and the aforementioned idea,

$$\int_V \exp(-\beta|\pmb{p}-\pmb{r}|)\frac{\cos\phi_p}{|\pmb{p}-\pmb{r}|^2}\chi(\pmb{p},\pmb{r})dV(\pmb{r}) = \int_S \left[\int_{L_{rp}} \exp(-\beta|\pmb{p}-\pmb{r}|)dL_{rp}(\pmb{r}')\right]\frac{\cos\phi_p \cos\phi_r}{|\pmb{p}-\pmb{r}|^2}\chi(\pmb{p},\pmb{r})dS(\pmb{r})$$

$$= \int_{S'}\left[\int_{L_{rp}} \exp(-\beta|\pmb{p}-\pmb{r}'|)dL_{rp}(\pmb{r}')\right]\frac{\cos\phi_p \cos\phi_r}{|\pmb{p}-\pmb{r}|^2}dS(\pmb{r})$$

$$= \frac{1}{\beta}\int_{S'}\left[1-\exp(-\beta|\pmb{p}-\pmb{r}|)\right]\frac{\cos\phi_p \cos\phi_r}{|\pmb{p}-\pmb{r}|^2}dS(\pmb{r})$$

$$\leq \frac{1}{\beta}\int_{S'} \frac{\cos\phi_p \cos\phi_r}{|\pmb{p}-\pmb{r}|^2}dS(\pmb{r}).$$

Adopting the above last inequality, we obtain

$$\int_V \exp(-\beta|\pmb{p}-\pmb{r}|)\frac{\cos\phi_p}{|\pmb{p}-\pmb{r}|^2}\chi(\pmb{p},\pmb{r})dV(\pmb{r}) \leq \frac{\pi}{\beta}.$$

The proofs are finished.

**Lemma 3**. Assume $S$ is a piecewise smooth surface of bounded domain $V$ in $R^3$. Let $\pmb{p} \in S$ be a point at which $S$ is smooth. Let $R$ denote the diameter of $V$. Then

$$\int_V \exp(-\beta|\pmb{p}-\pmb{r}|)\frac{1}{|\pmb{p}-\pmb{r}|^2}\chi(\pmb{p},\pmb{r})dV(\pmb{r}) = \frac{4\pi}{\beta}[1-\exp(-\beta R)], \qquad (11)$$

$$\int_{S} \exp(-\beta|\boldsymbol{p}-\boldsymbol{r}|) \frac{\cos\phi_r}{|\boldsymbol{p}-\boldsymbol{r}|^2} \chi(\boldsymbol{p},\boldsymbol{r}) dS(\boldsymbol{r}) \leq 4\pi. \tag{12}$$

hold.

*Proof.* As the preceding proofs in Lemma 2, we always can find a bounded domain $V' \subset V$ with boundary $S'$, in which $\boldsymbol{p}$ can see all the other point in the closure $\overline{V}'$, and

$$\int_{V} \exp(-\beta|\boldsymbol{p}-\boldsymbol{r}|) \frac{1}{|\boldsymbol{p}-\boldsymbol{r}|^2} \chi(\boldsymbol{p},\boldsymbol{r}) dV(\boldsymbol{r}) = \int_{V'} \exp(-\beta|\boldsymbol{p}-\boldsymbol{r}|) \frac{1}{|\boldsymbol{p}-\boldsymbol{r}|^2} dV(\boldsymbol{r}).$$

Let $\delta$ be a sufficiently small value and denote by $B(\boldsymbol{p};\delta)$ an $\delta$-neighborhood of $\boldsymbol{p}$. Exclude $B(\boldsymbol{p};\delta)$ from $V'$, and denote the remaining domain by $V'_\delta$.

Then, we have

$$\int_{V'} \exp(-\beta|\boldsymbol{p}-\boldsymbol{r}|) \frac{1}{|\boldsymbol{p}-\boldsymbol{r}|^2} dV(\boldsymbol{r}) = \int_{B(\boldsymbol{p};\delta)} \exp(-\beta|\boldsymbol{p}-\boldsymbol{r}|) \frac{1}{|\boldsymbol{p}-\boldsymbol{r}|^2} dV(\boldsymbol{r})$$

$$+ \int_{V'_\delta} \exp(-\beta|\boldsymbol{p}-\boldsymbol{r}|) \frac{1}{|\boldsymbol{p}-\boldsymbol{r}|^2} dV(\boldsymbol{r}).$$

For the first integral of the right hand side,

$$\int_{B(\boldsymbol{p};\delta)} \exp(-\beta|\boldsymbol{p}-\boldsymbol{r}|) \frac{1}{|\boldsymbol{p}-\boldsymbol{r}|^2} dV(\boldsymbol{r}) = \int_{B(\boldsymbol{0};\delta)} \exp(-\beta|\boldsymbol{r}|) \frac{1}{|\boldsymbol{r}|^2} dV(\boldsymbol{r})$$

$$= \int_0^\delta \left[ \int_{\partial B(\boldsymbol{0};r)} \exp(-\beta|\boldsymbol{r}|) \frac{1}{|\boldsymbol{r}|^2} dS \right] dr$$

$$= 4\pi \int_0^\delta \exp(-\beta r) dr$$

$$= \frac{4\pi}{\beta} \left[1 - \exp(-\beta\delta)\right].$$

Thus this integral goes to zero as $\delta \to 0$.

For the second integral,

$$\int_{V'_\delta} \exp(-\beta|\boldsymbol{p}-\boldsymbol{r}|) \frac{1}{|\boldsymbol{p}-\boldsymbol{r}|^2} dV(\boldsymbol{r}) \leq \int_{B(\boldsymbol{p};R)\setminus B(\boldsymbol{p};\delta)} \exp(-\beta|\boldsymbol{p}-\boldsymbol{r}|) \frac{1}{|\boldsymbol{p}-\boldsymbol{r}|^2} dV(\boldsymbol{r})$$

$$= \int_{B(\boldsymbol{0};R)\setminus B(\boldsymbol{0};\delta)} \exp(-\beta|\boldsymbol{r}|) \frac{1}{|\boldsymbol{r}|^2} dV(\boldsymbol{r})$$

$$= \int_\delta^R \left[ \int_{\partial B(\boldsymbol{0};r)} \exp(-\beta|\boldsymbol{r}|) \frac{1}{|\boldsymbol{r}|^2} dS \right] dr$$

$$= 4\pi \int_\delta^R \exp(-\beta r) dr$$

$$= \frac{4\pi}{\beta} \left[\exp(-\beta\delta) - \exp(-\beta R)\right].$$

So as $\delta \to 0$,

$$\lim_{\delta \to 0} \int_{V'_\delta} \exp(-\beta|\boldsymbol{p}-\boldsymbol{r}|) \frac{1}{|\boldsymbol{p}-\boldsymbol{r}|^2} dV(\boldsymbol{r}) \le \frac{4\pi}{\beta}\left[1-\exp(-\beta R)\right].$$

To sum up,

$$\int_V \exp(-\beta|\boldsymbol{p}-\boldsymbol{r}|) \frac{1}{|\boldsymbol{p}-\boldsymbol{r}|^2} \chi(\boldsymbol{p},\boldsymbol{r}) dV(\boldsymbol{r}) \le \frac{4\pi}{\beta}\left[1-\exp(-\beta R)\right].$$

Adopting the same idea to eliminate the Boolean function, that is,

$$\int_S \frac{\cos\phi_r}{|\boldsymbol{p}-\boldsymbol{r}|^2} \chi(\boldsymbol{p},\boldsymbol{r}) dS(\boldsymbol{r}) = \int_{S'} \frac{\cos\phi_r}{|\boldsymbol{p}-\boldsymbol{r}|^2} dS(\boldsymbol{r}).$$

Let $S'_\delta$ denote the boundary of $V'_\delta$, then

$$\int_{S'} \frac{\cos\phi_r}{|\boldsymbol{r}-\boldsymbol{p}|^2} dS(\boldsymbol{r}) = \int_{S'_\delta} \frac{(\boldsymbol{r}-\boldsymbol{p})\cdot\boldsymbol{n}_r}{|\boldsymbol{r}-\boldsymbol{p}|^3} dS(\boldsymbol{r}) + \int_{\partial B(\boldsymbol{p},\delta)} \frac{(\boldsymbol{r}-\boldsymbol{p})\cdot\boldsymbol{n}_r}{|\boldsymbol{r}-\boldsymbol{p}|^3} dS(\boldsymbol{r}).$$

Applying the divergence theorem to the first integral on the right hand side,

$$\int_{S'_\delta} \frac{(\boldsymbol{r}-\boldsymbol{p})\cdot\boldsymbol{n}_r}{|\boldsymbol{r}-\boldsymbol{p}|^3} dS(\boldsymbol{r}) = -\int_{V'_\delta} \nabla\cdot\left(\frac{\boldsymbol{r}-\boldsymbol{p}}{|\boldsymbol{r}-\boldsymbol{p}|^3}\right) dV(\boldsymbol{r}).$$

A straightforward computation shows

$$\nabla\cdot\left(\frac{\boldsymbol{r}-\boldsymbol{p}}{|\boldsymbol{r}-\boldsymbol{p}|^3}\right) = 0, \quad \text{for } \boldsymbol{r} \in V'_\delta.$$

And, therefore,

$$\int_{S'_\delta} \frac{(\boldsymbol{r}-\boldsymbol{p})\cdot\boldsymbol{n}_r}{|\boldsymbol{r}-\boldsymbol{p}|^3} dS(\boldsymbol{r}) = 0.$$

Besides, for $\boldsymbol{r} \in \partial B(\boldsymbol{p},\delta)$,

$$\boldsymbol{n}_r = \frac{\boldsymbol{r}-\boldsymbol{p}}{|\boldsymbol{r}-\boldsymbol{p}|}, \quad \boldsymbol{n}_r \cdot \frac{\boldsymbol{r}-\boldsymbol{p}}{|\boldsymbol{r}-\boldsymbol{p}|} = 1.$$

Then,

$$\int_{\partial B(\boldsymbol{p},\delta)} \frac{(\boldsymbol{r}-\boldsymbol{p})\cdot\boldsymbol{n}_r}{|\boldsymbol{r}-\boldsymbol{p}|^3} dS(\boldsymbol{r}) = \int_{\partial B(\boldsymbol{p},\delta)} \frac{1}{|\boldsymbol{r}-\boldsymbol{p}|^2} dS(\boldsymbol{r})$$

$$= \frac{1}{\delta^2} \int_{\partial B(\boldsymbol{p},\delta)} 1 dS(\boldsymbol{r})$$

$$= 4\pi.$$

So, as a consequence of above results,

$$\int_S \frac{\cos\phi_r}{|\boldsymbol{r}-\boldsymbol{p}|^2} dS(\boldsymbol{r}) = 4\pi.$$

That is,

$$\int_S \exp(-\beta|\boldsymbol{p}-\boldsymbol{r}|)\frac{\cos\phi_r}{|\boldsymbol{p}-\boldsymbol{r}|^2}\chi(\boldsymbol{p},\boldsymbol{r})dS(\boldsymbol{r}) \leq 4\pi.$$

The proof is finished.

**Corollary 1.** The operators $K_1, K_2, K_3, K_4$ are non-negative, and

$$\|K_1\|_{L^p(S)} < 1-\varepsilon, \quad \|K_2\|_{L^p(S)} < \frac{\varepsilon\sigma_s}{4\beta},$$

$$\|K_3\|_{L^p(V)} < \frac{\sigma_s}{\beta}[1-\exp(-\beta R)], \quad \|K_4\|_{L^p(V)} < \frac{4(1-\varepsilon)}{\varepsilon}.$$

for $p \in [1,\infty]$.

Proof. Let first

$$k_1(\boldsymbol{p},\boldsymbol{r}) = \exp[-\beta|\boldsymbol{p}-\boldsymbol{r}|]\frac{\cos\phi_p \cos\phi_r}{|\boldsymbol{p}-\boldsymbol{r}|^2}\chi(\boldsymbol{p},\boldsymbol{r}).$$

Suppose $1 < p < \infty$ and $u \in L^p(S)$. Then, following [30]

$$|K_1 u(\boldsymbol{p})| = \left|\int_S k_1^{1/p+1/q}(\boldsymbol{p},\boldsymbol{r})u(\boldsymbol{r})d(\boldsymbol{r})\right|$$

$$\leq \left(\int_S k_1(\boldsymbol{p},\boldsymbol{r})d(\boldsymbol{r})\right)^{1/q}\left(\int_S k_1(\boldsymbol{p},\boldsymbol{r})|u(\boldsymbol{r})|^p d(\boldsymbol{r})\right)^{1/p}$$

Then,

$$\int_S |K_1 u(\boldsymbol{p})|^p d\boldsymbol{p} \leq (1-\varepsilon)^{p/q}\int_S\int_S k_1(\boldsymbol{p},\boldsymbol{r})|u(\boldsymbol{r})|^p d\boldsymbol{r}d\boldsymbol{p}$$

$$= (1-\varepsilon)^{p/q}\int_S |u(\boldsymbol{r})|^p \int_S k_1(\boldsymbol{p},\boldsymbol{r})d\boldsymbol{p}d\boldsymbol{r}$$

$$= (1-\varepsilon)^{1+p/q}\int_S |u(\boldsymbol{r})|^p d\boldsymbol{r}$$

So, we have

$$\|K_1\|_{L^p(S)} = \sup_{\substack{u\in L^p(S)\\u\neq 0}}\frac{\|K_1 u\|_{L^p(S)}}{\|u\|_{L^p(S)}} \leq (1-\varepsilon)$$

The cases with $p=1$ and $p=\infty$ are straightforward.
So similarly, the other three inequalities are hold.

For actual physical problem, the emissivity satisfies
$$0 < \varepsilon < 1.$$
Then, we always have

$$\|K_1\|_{L^p(S)} < 1, \quad \|K_3\|_{L^p(V)} < 1. \tag{13}$$

With the help of above properties of the four operators, we are now in a position to prove the existence and uniqueness of solutions to Eqs. (4).

**Theorem 1.** Assume that there exists a constant $0 < \varepsilon_0 < 1$ such that $0 < \varepsilon_0 \leq \varepsilon \leq 1$. If $\dfrac{\sigma_s}{\beta + \sigma_s} < \varepsilon_0$, then Eqs. (4) has a unique solution.

*Proof.* The existence and uniqueness of solution to Eqs. (4) is equivalent to the reversibility of operator matrix

$$\begin{pmatrix} I - K_1 & -K_2 \\ -K_4 & I - K_3 \end{pmatrix}.$$

Using the inequalities (13) and the geometric series theorem, the operators $I - K_1$ and $I - K_3$ are invertible. Then

$$\begin{pmatrix} I - K_1 & -K_2 \\ -K_4 & I - K_3 \end{pmatrix} = \begin{pmatrix} I - K_1 & -K_2 \\ 0 & (I - K_3)\left[I - (I - K_3)^{-1} K_2 (I - K_1)^{-1} K_4\right] \end{pmatrix}.$$

Now, we only need to show the reversibility of the operator $I - (I - K_3)^{-1} K_2 (I - K_1)^{-1} K_4$.

Adopting the Neumann series theorem, we have

$$\left\|(I - K_1)^{-1}\right\| \leq \frac{1}{1 - \|K_1\|}, \quad \left\|(I - K_3)^{-1}\right\| \leq \frac{1}{1 - \|K_3\|}.$$

Then

$$\left\|(I - K_3)^{-1} K_2 (I - K_1)^{-1} K_4\right\| \leq \frac{\|K_2\|\|K_4\|}{(1 - \|K_1\|)(1 - \|K_3\|)}$$

$$\leq \frac{(1 - \varepsilon)}{\varepsilon} \frac{\sigma_s}{\sigma_a + \sigma_s \exp(-\beta R)}$$

$$\leq \frac{(1 - \varepsilon_0)}{\varepsilon_0} \frac{\sigma_s}{\sigma_a + \sigma_s \exp(-\beta R)}.$$

According to the Neumann series theorem, the operator $I - (I - K_3)^{-1} K_2 (I - K_1)^{-1} K_4$ is reversibility if $\left\|(I - K_3)^{-1} K_2 (I - K_1)^{-1} K_4\right\| < 1$, that is,

$$\frac{(1 - \varepsilon_0)}{\varepsilon_0} \frac{\sigma_s}{\sigma_a + \sigma_s \exp(-\beta R)} < 1$$

This equality is equivalent to

$$\frac{(1 - \varepsilon_0)}{\varepsilon_0} < \frac{\sigma_a}{\sigma_s} + \exp(-\beta R)$$

Further, we have

$$\frac{\sigma_a}{\sigma_s} + \exp(-\beta R) < \frac{\sigma_a}{\sigma_s} + 1 = \frac{\beta}{\sigma_s}$$

So, if

$$\frac{\sigma_s}{\beta + \sigma_s} < \varepsilon_0$$

holds, then the reversibility of the operator $I - (I - K_3)^{-1} K_2 (I - K_1)^{-1} K_4$ is obtained.

To sum up, we obtain that

$$\begin{pmatrix} I - K_1 & -K_2 \\ -K_4 & I - K_3 \end{pmatrix}$$

is invertible.

The proof is finished.

We note that the assumptions in Theorem 1 are satisfied in most situations. First, the theorem is valid for radiative heat transfer within black-body boundary or pure absorbing-emitting media. Besides, when the scattering coefficient is enough small, the theorem 1 is suitable for any situations with diffuse and gray boundary. For the situation with black boundary, the solvability of the RITEs in any semitransparent media is guaranteed by Theorem 1.

## 4. The iterative scheme and Numerical implementation

In order to simulate the coupled integral system involved in the radiative heat flux and the incident energy, an iterative scheme is constructed as

$$\begin{cases} q^n - K_1 q^n = K_2 G^n + f \\ G^{n+1} = K_3 G^n + K_4 q^n + g \end{cases} \quad (14)$$

Corollary 2. Under the hypotheses in Theorem 1, the iterative scheme (14) is convergent.

Proof. Solve the first equation in (14), we have

$$q^n = (I - K_1)^{-1} K_2 G^n + (I - K_1)^{-1} f,$$

where we take advantage of the invertibility of operator $I - K_1$. Substituting this formulation into the second equation in system (10), we have

$$G^{n+1} = K_3 G^n + K_4 \left[ (I - K_1)^{-1} K_2 G^n + (I - K_1)^{-1} f \right] + g.$$

Rearranging the right hand side,

$$G^{n+1} = \left[ K_3 + K_4 (I - K_1)^{-1} K_2 \right] G^n + \tilde{g}.$$

From the above equation, the error equation can be obtained as

$$e^{n+1} = \left[ K_3 + K_4 (I - K_1)^{-1} K_2 \right] e^n.$$

Then the convergence of iterative scheme (10) is equivalent to the contractibility of the operator $K_3 + K_4 (I - K_1)^{-1} K_2$.

$$\left\| K_3 + K_4 (I - K_1)^{-1} K_2 \right\| \leq \left\| K_3 \right\| + \left\| K_4 (I - K_1)^{-1} K_2 \right\|$$

$$\leq \frac{\sigma_s}{\beta} \left[ 1 - \exp(-\beta R) \right] + \frac{\sigma_s}{\beta} \frac{(1-\varepsilon)}{\varepsilon}$$

$$= \frac{\sigma_s}{\beta} \left[ \frac{1}{\varepsilon} - \frac{1}{\exp(\beta R)} \right]$$

$$< \frac{\sigma_s}{\beta} \frac{1-\varepsilon_0}{\varepsilon_0}.$$

From the assumption,

$$\left\| K_3 + K_4 (I - K_1)^{-1} K_2 \right\| < 1.$$

That is, the iteration scheme (10) is convergent.

The corollary 2 shows that if the solution of RITEs is unique, then the iteration scheme (15) always convergent.

BEM, as a technique of discretization of integral equations, is the proper tool to handle with RITEs. In this paper, BEM based on the collocation scheme is used to discretize system (1).

Firstly, substituting the geometric relationship [21]

$$dL_{rp}(\boldsymbol{r}') \frac{\cos \phi_r dS(\boldsymbol{r})}{|\boldsymbol{r} - \boldsymbol{p}|^2} = \frac{dV(\boldsymbol{r}')}{|\boldsymbol{r}' - \boldsymbol{p}|^2}$$

into Eqs. (1) (2), the alternative integral formulae of radiative transfer can be written as

$$q(\boldsymbol{p}) + \varepsilon E_b(\boldsymbol{p}) = \varepsilon \int_S \left[ E_b(\boldsymbol{r}) + \frac{1-\varepsilon}{\varepsilon} q(\boldsymbol{r}) \right] P_1(\boldsymbol{p},\boldsymbol{r}) dS(\boldsymbol{r})$$
$$+ \varepsilon \int_S \left[ \int_{L_{rp}} I_b(\boldsymbol{r}') \exp(-\beta |\boldsymbol{p} - \boldsymbol{r}'|) d\boldsymbol{r}' \right] P_2(\boldsymbol{p},\boldsymbol{r}) dS(\boldsymbol{r}) \qquad (15)$$
$$+ \varepsilon \int_S \left[ \int_{L_{rp}} G(\boldsymbol{r}') \exp(-\beta |\boldsymbol{p} - \boldsymbol{r}'|) d\boldsymbol{r}' \right] P_3(\boldsymbol{p},\boldsymbol{r}) dS(\boldsymbol{r}), \qquad \boldsymbol{p} \in S,$$

$$G(\boldsymbol{p}) = \int_S \left[ E_b(\boldsymbol{r}) + \frac{1-\varepsilon_r}{\varepsilon_r} q(\boldsymbol{r}) \right] P_4(\boldsymbol{p},\boldsymbol{r}) dS(\boldsymbol{r})$$
$$+ \int_S \left[ \int_{L_{rp}} I_b(\boldsymbol{r}') \exp(-\beta |\boldsymbol{p} - \boldsymbol{r}'|) d\boldsymbol{r}' \right] P_5(\boldsymbol{p},\boldsymbol{r}) dS(\boldsymbol{r}) \qquad (16)$$
$$+ \int_S \left[ \int_{L_{rp}} G(\boldsymbol{r}') \exp(-\beta |\boldsymbol{p} - \boldsymbol{r}'|) d\boldsymbol{r}' \right] P_6(\boldsymbol{p},\boldsymbol{r}) dS(\boldsymbol{r}), \qquad \boldsymbol{p} \in V,$$

where the kernel functions are defined as

$$P_1(\boldsymbol{p},\boldsymbol{r}) = \exp(-\beta|\boldsymbol{p}-\boldsymbol{r}|)\frac{\cos\phi_p \cos\phi_r}{\pi|\boldsymbol{p}-\boldsymbol{r}|^2}\chi(\boldsymbol{p},\boldsymbol{r})$$

$$P_2(\boldsymbol{p},\boldsymbol{r}) = \sigma_a \frac{\cos\phi_p \cos\phi_r}{|\boldsymbol{p}-\boldsymbol{r}|^2}\chi(\boldsymbol{p},\boldsymbol{r})$$

$$P_3(\boldsymbol{p},\boldsymbol{r}) = \frac{\sigma_s}{4\pi}\frac{\cos\phi_p \cos\phi_r}{|\boldsymbol{p}-\boldsymbol{r}|^2}\chi(\boldsymbol{p},\boldsymbol{r}) \qquad (17)$$

$$P_4(\boldsymbol{p},\boldsymbol{r}) = \exp(-\beta|\boldsymbol{p}-\boldsymbol{r}|)\frac{\cos\phi_r}{\pi|\boldsymbol{p}-\boldsymbol{r}|^2}\chi(\boldsymbol{p},\boldsymbol{r})$$

$$P_5(\boldsymbol{p},\boldsymbol{r}) = \sigma_a \frac{\cos\phi_r}{|\boldsymbol{p}-\boldsymbol{r}|^2}\chi(\boldsymbol{p},\boldsymbol{r})$$

$$P_6(\boldsymbol{p},\boldsymbol{r}) = \frac{\sigma_s}{4\pi}\frac{\cos\phi_r}{|\boldsymbol{p}-\boldsymbol{r}|^2}\chi(\boldsymbol{p},\boldsymbol{r})$$

Then, the discretization of the boundary integral equations (13) (14) is carried out by using the standard boundary element discretization technique. First, the boundary bounding the domain under consideration is discretized to a set of plane elements denoted by $\{\Delta S_k\}_{k=1}^{N_e}$. And approximate $q$ by linear interpolation on each element. For example, on element $\Delta S_k$, we have

$$q(\boldsymbol{p}) = q(\xi,\eta) = \sum_{\alpha=1}^{M} F_\alpha(\xi,\eta) q^\alpha$$

Where $\xi,\eta$ are the intrinsic coordinate, and $F_\alpha(\xi,\eta)$, $q^\alpha$ denote the bilinear shape function and the nodal value of $q$ at local nodes $\alpha$ on the element $\Delta S_k$. Denote by $M$ the number of nodal on each element.

Let $\boldsymbol{p}$ in Eq. (13) to be given (in turn) by $\boldsymbol{p}_i$ for $i=1,\cdots,N_p$ (where $N_p$ is the number of collocation points on the boundary). Let $\boldsymbol{p}$ in Eq. (14) to be given (in turn) by $\boldsymbol{p}_j$ for $j=1,\cdots,N_i$ (where $N_i$ is the number of collocation points in the interior of the domain). We can obtain the approximate system

$$q_i - \sum_{k=1}^{Nel}\sum_{\alpha=1}^{M} g_{ik}^\alpha q^\alpha = \sum_{k=1}^{Nel}\sum_{l=1}^{L} f_{ik}^l S^l + h_i, \qquad i=1,2,\cdots,N_p, \qquad (18)$$

$$G_j = \sum_{k=1}^{Nel}\sum_{l=1}^{L} u_{jk}^l G^l + \sum_{k=1}^{Nel}\sum_{\alpha=1}^{M} v_{jk}^\alpha q^\alpha + t_j, \qquad j=1,2,\cdots,N_i, \qquad (19)$$

where

$$g_{ik}^{\alpha} = (1-\varepsilon)\int_{\Delta S_k} F_\alpha P_1(\mathbf{p}_i,\mathbf{r})dS(\mathbf{r}),$$

$$f_{ik}^{l} = \varepsilon\int_{\Delta S_k}\left[\int_{L_{rp}^l}\exp(-\beta|\mathbf{p}_i-\mathbf{r}'|)d\mathbf{r}'\right]P_3(\mathbf{p}_i,\mathbf{r})dS(\mathbf{r}),$$

$$h_i = \varepsilon\sum_{k=1}^{Nel}\int_{\Delta S_k}\sum_{\alpha=1}^{M}F_\alpha E_b^\alpha(\mathbf{r})P_1(\mathbf{p}_i,\mathbf{r})dS(\mathbf{r}) - \varepsilon E_b(\mathbf{p}_i) +$$

$$\varepsilon\sum_{k=1}^{Nel}\int_{\Delta S_k}\sum_{l=1}^{L}\left[\int_{L_{rp}^l}I_b(\mathbf{r}')\exp(-\beta|\mathbf{p}_i-\mathbf{r}'|)d\mathbf{r}'\right]P_2(\mathbf{p},\mathbf{r})dS(\mathbf{r}),$$

$$u_{jk}^{l} = \int_{\Delta S_k}\left[\int_{L_{rp}^l}\exp(-\beta|\mathbf{p}_j-\mathbf{r}'|)d\mathbf{r}'\right]P_6(\mathbf{p}_j,\mathbf{r})dS(\mathbf{r}),$$

$$v_{jk}^{\alpha} = \left(\frac{1-\varepsilon}{\varepsilon}\right)\int_{\Delta S_k}F_\alpha P_4(\mathbf{p}_j,\mathbf{r})dS(\mathbf{r}),$$

$$t_j = \sum_{k=1}^{Nel}\int_{\Delta S_k}\sum_{\alpha=1}^{M}F_\alpha E_b^\alpha(\mathbf{r})P_4(\mathbf{p},\mathbf{r})dS(\mathbf{r}) +$$

$$\sum_{k=1}^{Nel}\int_{\Delta S_k}\sum_{l=1}^{L}\left[\int_{L_{rp}^l}I_b(\mathbf{r}')\exp(-\beta|\mathbf{p}_j-\mathbf{r}'|)d\mathbf{r}'\right]P_5(\mathbf{p}_j,\mathbf{r})dS(\mathbf{r}). \quad (20)$$

$q_i, G_j$ are the abbreviations of $q(\mathbf{p}_i), G(\mathbf{p}_j)$. $\{L_{rp}^l\}_{l=1}^{L}$ denotes the subdivision of the line $L_{rp}$.

Then system (18)-(19) are solved by the iterative scheme (14).

Additionally, the computational domain is covered by a regular grid of cuboid cells in order to calculate the line integral. When the temperature and incidence radiation within each small cuboid cell are assumed to be constants, the line integrals in $f_{ik}^l, h_i, u_{jk}^l, t_j$ can be evaluated analytically.

The nearly singular term, $\dfrac{1}{|\mathbf{p}-\mathbf{r}|^2}$, can be handled by the technique proposed by Eberwien et al. [26]. The idea allows us to ensure that the error made by numerical integration is nearly constant, regardless of the proximity of the source point. The detailed procedure is not repeated here. The interests can be referred to [26].

Note too that the Boolean function in kernels (15) needs to be determined by an effectively visibility algorithm. This procedure is very time-consuming. In addition, the accurate detection of the limits of the zones is a hercules task and is crucial for the accuracy of the result of the numerical simulation. The aim of next section is to develop a high-precision visibility algorithm.

**5. The algorithm for the shadow detection**

The present algorithm can be widely used to solve 2D and 3D problems. Moreover, the subdivision of the current element bases on the detection of obstructing nodes. In the limit, this algorithm can capture the shadow domain exactly.

Now the detailed detection and subdivision procedure is described as following:

Note that the developed algorithm is valid for flat elements. Denote by $\mathbf{p}$ the source point. $S_k, S_b$ represent the current element (potential active element) and the third party element

(potential obstructing element), respectively. Meanwhile, when $p$ locates at the boundary, let $n_p, n_k$ denote the unit outward normal vector of $p, S_k$, respectively. Let $r, r_1, r_2$ denote the middle point and the two endpoints of $S_k$.

**Step 1. Creation of the list of active elements**

This list is updated as the procedure proceeds.

First, examine the signs of

$$d_1 = n_p \cdot (p - r), \ d_2 = n_k \cdot (r - p)$$

If the signs of $d_1, d_2$ are both positive, the current element $S_k$ is considered as to be a potential active element and is added to the list of active elements. Repeating this operator for all boundary elements, we build the list of active elements for $p$.

**Step 2. Creation of the list of blocking element**

For any a potential active element, the list of blocking elements is established by scanning all the other elements with the following operations.

First, the elements with the source point or coplanar potential active element are excluded. After that, a cylinder window is created as Fig. 1, of which the diameter is the sum of diameters of potential active and potential blocking element, and the height is $|p - r|$.

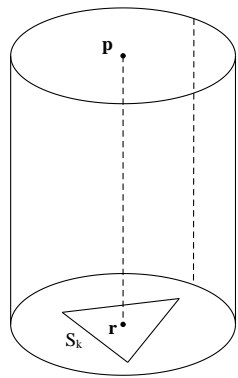

Figure. 1.

Judge whether the center of $S_b$ lies in this cylinder or not. If so, $S_b$ is needed to do further judgment. If not, the judgment of this potential blocking element is terminated.

Then, a cone window is constructed as Fig. 2. The bottom diameter of the cone is equal to the diameter of $S_k$ and its height is $|p - r|$.

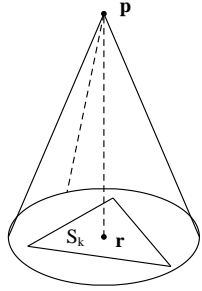

Figure. 2.

Judge whether one of the vertexes of $S_j$ lies in this cone or not. If so, put this third party element $S_j$ into the list of blocking elements about $p$ and $S_k$. If not, more judgment is needed to exclude the situation as shown in Fig. 3, which are described as follows.

First, judge whether the line segment $pr$ intersects $S_j$ or not. If so, we conclude that $S_k$ cannot be viewed from $p$ immediately. $S_k$ is discarded from the list of potential active elements and to check the next element. If not, judge whether one of the lines $pr_1, pr_2, pr_3$ intersects $S_j$ or not. If it is true, put $S_j$ into the list of blocking elements. If not, $S_j$ is discarded and the next third party element is continuous.

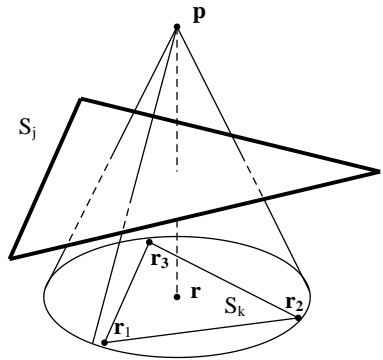

Figure. 3.

A special case is Sk may be shaded by the union of some adjacent blocking element. To avoid this misjudgement, we use a simple operation, which is to judge whether all the lines $pr_1, pr_2, pr_3$ intersect some blocking element in the list of blocking elements. If so, this active element $S_k$ is discarded from the list of active elements. If not, no operation is performed.

**Step 3. Subdivision of the active element**

By previous steps, a list of active elements for point $p$ is established. Meanwhile, a list of blocking elements about $p$ and any active element is established. For one active element, if its list of blocking elements is empty, the numerical integration is performed directly. For the remaining active element, this step is essential, that is, a subdivision is proceeds as follows.

First, this element is subdivided into four sub-elements as shown in Fig. 4.

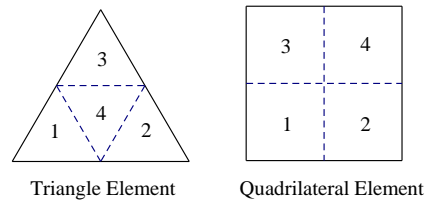

Figure. 4. Subdivision of element.

For all the sub-elements, the step 2 above is performed running through all blocking elements. When some sub-element is blocked partly, the subdivision is continuous. For sub-element fully visible, the integration on this sub-element is performed. On the contrary, the sub-elements fully shaded are discarded.

The subdivision procedure is terminated when the area of the produced sub-element reaches a preset minimum value.

Theoretically, the above detection algorithm is exact when the preset minimum value takes the limit zero. Although the present algorithm is described based on 3D geometry, it can be applied to 2D geometries naturally. For 2D problems, the two windows in step 2 are rectangle and triangle, respectively.

## 6. Numerical Examples

Example 1. *L-shaped enclosure*

In this section, a 3D L-shaped enclosure is considered (see Fig. 5). The dimensions of this tested geometry are $W \times L \times H = 1 \times 3 \times 3 \left(m^3\right)$ (where $H_1 = 2m$, $L_1 = 1m$). The enclosure contains an emitting-absorbing medium at a temperature of 1000K. The walls are black at 500 K. Fig. 6 shows the effect of absorbing coefficient of the medium on the predicted net radiative heat flux along the *AA* line.

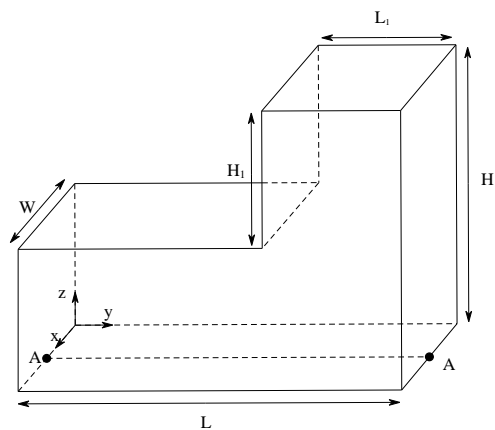

Figure. 5. The L-shaped enclosure.

The present model uses 2200 boundary discrete elements. The comparison with the existing results is shown in Fig. 6. From the results, the validity of present algorithm is verified. The proposed numerical procedure seems to be an attractive way of dealing with complex enclosure.

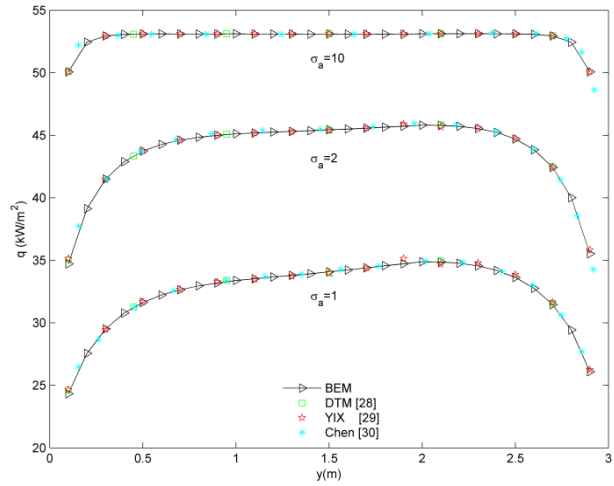

Figure. 6. Net radiative heat flux along line $AA$.

Example 2. *Cubic enclosures*

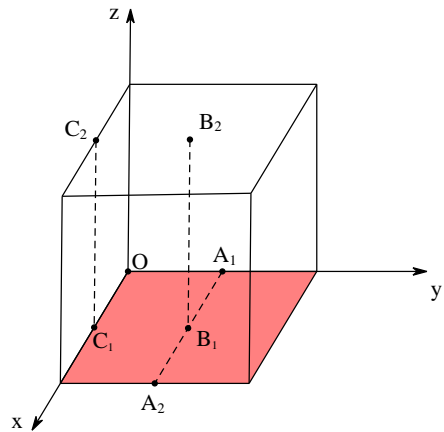

Figure. 7. System geometry.

The analysis of radiative transfer of multidimensional rectangular geometries has historically been associated with the design of combustion chambers and furnaces. This example refers to [26]. A cubic enclosure has unit length sides. All six walls of the cube are black. The scattering albedo of medium is equal to unity ($\Omega = \sigma_s / \beta = 1$). Like as [31], we also consider the dimensionless heat flux at $A_1A_2, C_1C_2$ and the average incident radiation at $B_1B_2$. The numerical comparison is shown in Fig. 8-10.

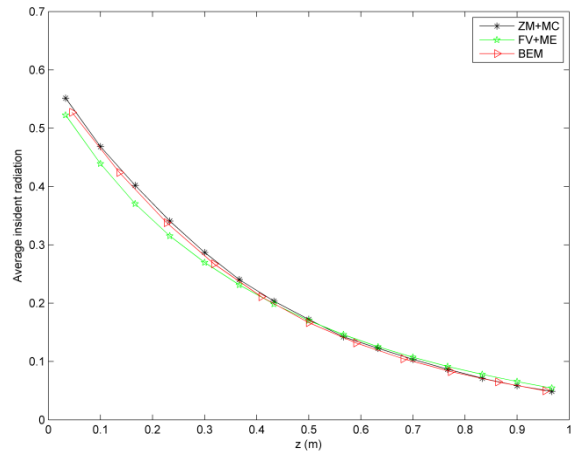

Figure. 8. Non-dimensional net radiative heat flux in z-direction along with $C_1C_2$.

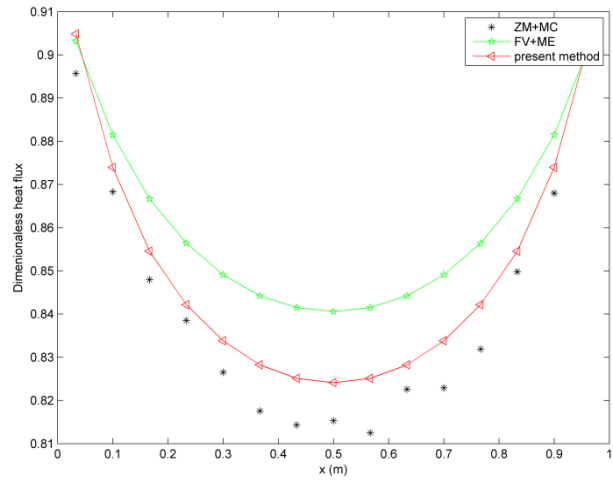

Figure. 9. Non-dimensional net radiative heat flux in x-direction along with $A_1A_2$.

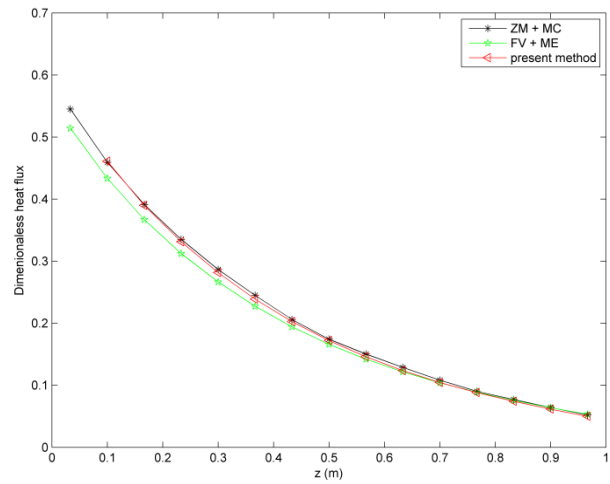

Figure. 10. Non-dimensional net radiative heat flux in z-direction along with $B_1B_2$.

The present model uses 1350 boundary elements and 1331 cuboid cells in the medium. From the comparative results, the present results are more identical to the references compared with the FVM.

**Conclusions**

Firstly, the properties of the four integral operators in RITEs are analyzed. Based on these results, the existence and uniqueness of solution to RITEs is proved and an iterative scheme, which is suggested for RITEs, is proved to be convergent. An improved high-precision visibility algorithm is developed to handle with non-convex domain. We also verify the effectiveness of scheme and the accuracy of algorithm by numerical examples. In consideration heavily used search operations in our algorithm, we are prepared to go further to improve efficiency of the method, which will be our future interest.

**Reference:**

[1] Agoshkov V. Boundary value problems for transport equations [M]. Springer Science & Business Media, 2012.

[2] Modest M F. Radiative heat transfer [M]. Academic press, 2013.

[3] Howell J R, Menguc M P, Siegel R. Thermal radiation heat transfer [M]. CRC press, 2010.

[4] Adams M L, Larsen E W. Fast iterative methods for discrete-ordinates particle transport calculations [J]. Progress in nuclear energy, 2002, 40(1): 3-159.

[5] Crosbie, A.L. and R.G. Schrenker, Exact expressions for radiative transfer in a three-dimensioanl rectangular geometry. Journal of Quantitative Spectroscopy & Radiative Transfer, 1982. 28(6): 507-526.

[6] Crosbie A L, Schrenker R G. Radiative transfer in a two-dimensional rectangular medium exposed to diffuse radiation [J]. Journal of Quantitative Spectroscopy and Radiative Transfer, 1984, 31(4): 339-372.

[7] Thynell ST. The integral form of the equation of transfer in finite, two-dimensional, cylindrical media. Journal of Quantitative Spectroscopy and Radiative Transfer, 1989; 42: 117–36.

[8] Tan Z. Radiative heat transfer in multidimensional emitting, absorbing, and anisotropic scattering media—mathematical formulation and numerical method [J]. Journal of heat transfer, 1989, 111(1): 141-147.

[9] Sun B, Zheng D, Klimpke B, et al. Modified boundary element method for radiative heat transfer analyses in emitting, absorbing and scattering media[J]. Engineering analysis with boundary elements, 1998, 21(2): 93-104.

[10] Li B Q, Cui X, Song S P. The Galerkin boundary element solution for thermal radiation problems [J]. Engineering Analysis with boundary elements, 2004, 28(7): 881-892.

[11] Altaç Z, Tekkalmaz M. Benchmark solutions of radiative transfer equation for three-dimensional rectangular homogeneous media [J]. Journal of Quantitative Spectroscopy and Radiative Transfer, 2008, 109(4): 587-607.

[12] Altac Z, Tekkalmaz M. Exact solution of radiative transfer equation for three-dimensional rectangular, linearly scattering medium [J]. Journal of Thermophysics and Heat Transfer, 2011, 25(2): 228-238.

[13] Atkinson K, Chandler G. The collocation method for solving the radiosity equation for unoccluded surfaces [J]. Journal of Integral Equations and Applications, 1998, 10(3): 253-290.

[14] Atkinson K, Chien D D K, Seol J. Numerical analysis of the radiosity equation using the collocation method [J]. Electronic Transactions on Numerical Analysis, 2000, 11: 94-120.

[15] Hansen O. The Local Behavior of the Solution of the Radiosity Equation at the Vertices of Polyhedral Domains in $R^3$ [J]. SIAM journal on mathematical analysis, 2001, 33(3): 718-750.

[16] Qatanani N A, Daraghmeh A. Asymptotic error analysis for the heat radiation boundary integral equation [J]. European

Journal of Mathematical Sciences, 2013, 2(1): 51-61.

[17] Tiihonen T. Stefan-Boltzmann radiation on non-convex surfaces [J]. Mathematical methods in the applied sciences, 1997, 20(1): 47-57.

[18] Laitinen M T, Tiihonen T. Integro-differential equation modelling heat transfer in conducting, radiating and semitransparent materials [J]. Mathematical methods in the applied sciences, 1998, 21(5): 375-392.

[19] Watt A. Fundamentals of three-dimensional computer graphics [M]. Addison-Wesley, 1989.

[20] Cohen M F, Wallace J R. Radiosity and realistic image synthesis [M]. Elsevier, 2012.

[21] Emery A F, Johansson O, Lobo M, et al. A comparative study of methods for computing the diffuse radiation view factors for complex structures [J]. Journal of heat transfer, 1991, 113(2): 413-422.

[22] Blobner J, Bialecki R A, Kuhn G. Boundary-element solution of coupled heat conduction-radiation problems in the presence of shadow zones [J]. Numerical Heat Transfer Part B Fundamentals, 2001, 39(5): 451-478.

[23] Li B Q, Cui X, Song S P. The Galerkin boundary element solution for thermal radiation problems [J]. Engineering Analysis with boundary elements, 2004, 28(7): 881-892.

[24] Kress R. Linear integral equations [M]. Springer, 2014.

[25] Viskanta R. Radiation transfer and interaction of convection with radiation heat transfer [J]. Advances in Heat Transfer, 1966, 3: 175-251.

[26] Eberwien U, Duenser C, Moser W. Efficient calculation of internal results in 2D elasticity BEM [J]. Engineering analysis with boundary elements, 2005, 29(5): 447-453.

[27] Białecki R A, Grela Ł. Application of the boundary element method in radiation [J]. Journal of Theoretical and Applied Mechanics, 1998, 36(2): 351-364.

[28] W. M. Malalasekera and E. H. James, Radiative Heat Transfer Calculations in Three-dimensional Complex Geometries, ASME J. Heat Transfer, vol. 118, no. 1, pp. 225–228, 1996.

[29] Hsu P, Tan Z. Radiative and combined-mode heat transfer within L-shaped nonhomogeneous and nongray participating media. Numerical Heat Transfer, Part A: Applications, 1997, 31(8): 819-835.

[30] Chen S S, Li B W, Tian X Y. Chebyshev collocation spectral domain decomposition method for coupled conductive and radiative heat transfer in a 3D L-shaped enclosure. Numerical Heat Transfer, Part B: Fundamentals, 2016, 70(3): 215-232.

[31] Trivic D N, Amon C H. Modeling the 3-D radiation of anisotropically scattering media by two different numerical methods [J]. International Journal of Heat and Mass Transfer, 2008, 51(11): 2711-2732.